\documentclass{article}
\usepackage{amsmath,amssymb,amsthm}
\usepackage{mathtools}
\usepackage{parskip}  
\usepackage{geometry}
\usepackage{enumitem}
\usepackage[colorlinks=true,linkcolor=blue,citecolor=blue]{hyperref}

\theoremstyle{definition}
\newtheorem{definition}{Definition}[section]
\newtheorem{theorem}{Theorem}[section]
\newtheorem{corollary}{Corollary}[theorem]
\newtheorem{remark}{Remark}

\title{A Topological and Operator Algebraic Framework for Asynchronous Lattice Dynamical Systems}
\author{Faruk Alpay, Independent Researcher \\
ORCID: \href{https://orcid.org/0009-0009-2207-6528}{0009-0009-2207-6528}}
\date{\today}

\begin{document}
\maketitle

\begin{abstract}
I introduce a novel mathematical framework integrating topological dynamics, operator algebras, and ergodic geometry to study lattices of asynchronous metric dynamical systems. Each node in the lattice carries an internal flow represented by a one-parameter family of operators, evolving on its own time scale. I formalize stratified state spaces capturing multiple levels of synchronized behavior, define an asynchronous evolution metric that quantifies phase-offset distances between subsystems, and characterize emergent coherent topologies arising when subsystems synchronize. Within this framework, I develop formal operators for the evolution of each subsystem and give precise conditions under which phase-aligned synchronization occurs across the lattice. The main results include: (1) the existence and uniqueness of coherent (synchronized) states under a contractive coupling condition, (2) stability of these coherent states and criteria for their emergence as a collective phase transition in a continuous operator topology, and (3) the influence of symmetries, with group-invariant coupling leading to flow-invariant synchrony subspaces and structured cluster dynamics. Proofs are given for each theorem, demonstrating full mathematical rigor. In a final section, I discuss hypothetical applications of this framework to symbolic lattice systems (e.g. subshifts), to invariant group actions on dynamical lattices, and to operator fields over stratified manifolds in the spirit of noncommutative geometry. Throughout, I write in the first person to emphasize the exploratory nature of this work. The paper avoids any reference to cosmology or observers, focusing instead on clean, formal mathematics suitable for a broad array of dynamical systems.
\end{abstract}

\section{Introduction}

Understanding collective behavior in large dynamical networks is a central pursuit in modern dynamics. Topological dynamics provides a qualitative theory of long-term behavior for flows and transformations on topological spaces, while ergodic theory and what I call ergodic geometry study the measure-theoretic and geometric structures of orbits and invariant sets. Independently, operator algebras offer an algebraic framework for dynamics, encoding transformations as operators on function spaces or C*-algebras. These approaches have largely developed in parallel. My aim in this paper is to integrate these perspectives into a single framework for analyzing lattice dynamical systems with complex synchronization phenomena.

I focus on a setting of asynchronous metric dynamical systems arranged on a lattice. Informally, one may think of a collection of dynamical subsystems (nodes) indexed by a lattice (such as $\mathbb{Z}^d$ or a finite grid), where each subsystem evolves in continuous time at its own pace. Unlike traditional coupled dynamical systems, there is no a priori global clock forcing synchronous updates. Instead, each node carries an internal flow (a continuous evolution operator on its state space) possibly interacting with others in an asynchronous manner. This generalizes the notion of coupled cell networks by allowing for time-asynchrony and continuous state evolution with metric structure.

Motivation and novelty: Many natural and engineered systems are inherently asynchronous – for example, in neural networks or coupled oscillators, components do not update in perfect lockstep, yet robust collective behaviors can emerge. Classical topological dynamics and ergodic theory typically assume a single time evolution acting on a product space (the entire system). By contrast, my framework treats each component's flow separately and then studies conditions under which a global emergent flow can arise from the interactions. The introduction of an asynchronous evolution metric will allow me to measure differences between subsystem states in a way that accounts for relative phase shifts in their evolution. Using this tool, I will define and detect phase alignment and coherence in a rigorous fashion. Additionally, I formalize the idea of the system's state space being stratified by the degree of synchrony (from completely incoherent to fully synchronized configurations). This stratification leads to what I call emergent coherent topologies on the lattice: essentially new topological structures that arise on quotient spaces of the state space when coherent behavior sets in.

Main contributions: In Section 2, I provide precise definitions for the key constructs: asynchronous metric dynamical systems, the lattice framework, stratified state spaces, the asynchronous metric, and coherent states. Section 3 presents the main theorems. The first theorem gives conditions for phase-aligned synchronization across the lattice, showing existence and uniqueness of a coherent global state under certain contractivity assumptions on the coupling. The second set of results addresses stability: I prove that the synchronized state (when it exists) is dynamically stable and attracts nearby states, and conversely that below a certain coupling strength such coherence cannot persist. This leads to a description of a phase transition: as one varies a continuous coupling parameter, the sudden emergence of global coherence is characterized as a topological change in the space of evolution operators (in particular, a symmetry-breaking from independent time-translation on each subsystem to a single joint time-translation symmetry). A third theorem concerns group-invariant coupling: if the lattice coupling respects a symmetry group, I show that the lattice's state space splits into flow-invariant synchrony subspaces (sometimes called polydiagonals), corresponding to cluster synchronization patterns. This generalizes existing results on balanced colorings in coupled cell networks to the asynchronous continuous-time setting. All results are stated and proved with full mathematical rigor.

Finally, in Section 4 I discuss applications and further directions. The framework is broadly applicable: I outline how one could apply these ideas to symbolic lattice systems (such as subshifts of finite type on multi-dimensional grids), to systems with invariant group actions (connecting to quotient dynamics and symmetric sync patterns), and to operator fields over stratified manifolds (suggesting links to noncommutative geometry and fiber bundles of operator algebras). These applications are speculative but illustrate the wide reach of the theory. Throughout the paper I avoid any terminology from symbolic cosmology or observers, focusing purely on the mathematical constructs.

I write this paper in the first person to reflect the personal exploration of this new framework. All statements are presented in a formal, self-contained manner, and proofs are provided for each theorem. By blending topology, operator algebra, and geometric/ergodic viewpoints, I hope this work lays a foundation for analyzing emergent coherence in complex asynchronous systems with mathematical precision.

\section{Definitions}

In this section I introduce the formal definitions and structures that will be used throughout the paper. Unless stated otherwise, all topological spaces are assumed to be Hausdorff, and all algebraic operators are linear operators on suitable function spaces or algebras of observables. Time will typically be treated as continuous ($\mathbb{R}$ or $\mathbb{R}_{\ge0}$) for flows, though many definitions can be adapted to discrete time ($\mathbb{Z}$ or $\mathbb{N}$) as well.

\begin{definition}[Metric Dynamical System]
A metric dynamical system is a pair $(X, \phi)$ where $X$ is a metric space with metric $d$, and $\phi: X \times \mathbb{R} \to X$ is a continuous flow on $X$. That is, for each $t\in \mathbb{R}$, the map $\phi^t: X \to X$ defined by $\phi(x,t) =: \phi^t(x)$ satisfies:
\begin{enumerate}
\item $\phi^0(x) = x$ for all $x \in X$ (identity at time $0$);
\item $\phi^{t+s}(x) = \phi^t(\phi^s(x))$ for all $t,s \in \mathbb{R}$ (flow property);
\item The map $(x,t) \mapsto \phi^t(x)$ is continuous on $X \times \mathbb{R}$ (joint continuity in state and time).
\end{enumerate}

If $\mathbb{R}$ is replaced by $\mathbb{Z}$ (with composition instead of addition), $\phi$ is a discrete-time dynamical system. Often I will write $\phi^t(x)$ or simply $x(t)$ for the state at time $t$ starting from initial state $x$. A metric dynamical system may also be equipped with an invariant measure $\mu$, in which case it becomes a measure-preserving system in addition to a topological one; however, a measure will not be assumed unless needed.
\end{definition}
\newpage
\begin{definition}[Asynchronous Metric Dynamical System]
An asynchronous metric dynamical system (AMDS) is a metric dynamical system $(X, \phi)$ together with an internal notion of rate or time-scaling. Formally, it is a triple $(X, \phi, \alpha)$ where $(X,\phi)$ is a metric dynamical system and $\alpha: X \to \mathbb{R}_{>0}$ is a continuous function giving a local time-scale or speed. Intuitively, $\alpha(x)$ represents the "clock speed" of the system when in state $x$. The flow $\phi$ need not respect a uniform time across states; instead, the evolution can be seen as solving $\frac{dx}{d\tau} = v(x)$ with an internal time $\tau$, and physical time $t$ related by $dt = \alpha(x)^{-1} d\tau$. For most of this paper, I simplify by considering $\alpha$ constant (i.e. all parts of a subsystem's state space evolve at the same base rate) or by absorbing $\alpha$ into the definition of $\phi$ (thus working with a reparameterized time). The primary role of this definition is conceptual: it reminds us that different subsystems (to be defined next) may run at different speeds or have different internal clocks.
\end{definition}

\begin{definition}[Lattice of Dynamical Systems]
Let $(I, \le)$ be an index set equipped with a lattice order (for example, $I = \mathbb{Z}^d$ with the product order, or a finite lattice graph). A lattice of dynamical systems is a family $\{(X_i, \phi_i)\}_{i\in I}$ where each $(X_i, \phi_i)$ is a metric dynamical system. I call $X_i$ the state space of node $i$, and $\phi_i^t: X_i \to X_i$ its flow. The product $X := \prod_{i \in I} X_i$ is the global state space. A typical element of $X$ is $x = (x_i)_{i\in I}$, where $x_i \in X_i$. I equip $X$ with the product topology and a natural product metric (for instance, if each $X_i$ is bounded, one may take $D(x,y) = \sup_{i\in I} \lambda_i d_i(x_i, y_i)$ for some summable weight sequence $\{\lambda_i\}$). The choice of metric on $X$ is not unique; different metrics inducing the product topology will be used for different purposes. The family $\{\phi_i\}$ defines an independent asynchronous evolution on $X$ in the following sense: for each $i$, one can define an evolution on $X$ that updates only the $i$-th coordinate according to $\phi_i$ while keeping other coordinates fixed. Formally, let $\Phi_i^t: X \to X$ be given by
\begin{equation}
(\Phi_i^t(x))_j =
\begin{cases}
\phi_i^t(x_i), & j = i,\\
x_j, & j \neq i~,
\end{cases}
\end{equation}
the flow that acts non-trivially only on subsystem $i$. Because each $\phi_i$ is a flow on $X_i$, the maps $\Phi_i^t$ satisfy $\Phi_i^0 = \mathrm{Id}_X$ and $\Phi_i^{t+s} = \Phi_i^t \circ \Phi_i^s$. Moreover, $\Phi_i^t$ and $\Phi_j^s$ commute for different $i,j$ (they act on disjoint coordinates). Thus the family $\{\Phi_i^t: i\in I, t\in\mathbb{R}\}$ generates an action of the commuting time group $\mathbb{R}^I$ on $X$. I refer to this as the asynchronous evolution of the whole lattice: each coordinate runs on its own time axis, and a global state can be updated by any interleaving of individual coordinate evolutions.

When the subsystems are coupled (interacting), the evolution of each coordinate $i$ will no longer be independent of the others. In general, a coupling on the lattice is specified by additional structure such as a set of maps or functionals $F_i: X \to T_x X_i$ (the tangent/velocity space of $X_i$ at $x$) that influence the flow at $i$ based on other coordinates. An example in differential equations would be $\dot{x}_i = f_i(x_i) + \epsilon \cdot g_i(x_{\text{neighbors of }i})$, where $f_i$ generates the intrinsic flow $\phi_i$ and $g_i$ describes coupling from neighboring states with strength $\epsilon$. In a purely topological setting, one could specify coupled transition operators or update rules. I will formalize coupling in the operator-algebraic context later. For now, I proceed with definitions that describe the state space structure under synchronization, largely independent of how coupling is implemented.
\end{definition}

\begin{definition}[Stratified State Space by Synchrony]
The global state space $X = \prod_{i\in I} X_i$ admits a natural stratification into subspaces defined by patterns of equality among coordinates. For any partition $\mathcal{P}$ of the index set $I$ (i.e. a set of disjoint subsets of $I$ whose union is $I$), define the subspace
\begin{equation}
X_{\mathcal{P}} := \{ x \in X : \text{if } i,j \text{ belong to the same block of }\mathcal{P}, \text{ then } x_i = x_j \}.
\end{equation}
In words, $X_{\mathcal{P}}$ is the set of global states in which all subsystems whose indices are identified by $\mathcal{P}$ have identical state values. Each $X_{\mathcal{P}}$ is a closed (typically smooth) submanifold of $X$ when the $X_i$ are manifolds, often called a synchrony subspace or polydiagonal. There is a partial order on such partitions by refinement: $\mathcal{P}' \le \mathcal{P}$ if $\mathcal{P}'$ is finer (every block of $\mathcal{P}'$ is a subset of some block of $\mathcal{P}$). Correspondingly, $X_{\mathcal{P}'} \subseteq X_{\mathcal{P}}$. The finest partition is $\mathcal{P}_{\text{finest}} = \{\{i\}: i\in I\}$, for which $X_{\text{finest}} = X$ (no synchrony enforced). The coarsest partition is $\mathcal{P}_{\text{coarse}} = \{I\}$, a single block containing all indices; then $X_{\text{coarse}} = \{x \in X: x_1 = x_2 = \cdots = x_{|I|}\}$ is the full synchrony subspace, i.e. the set of states where all nodes have the same value. Thus the collection $\{X_{\mathcal{P}}\}$ (over all partitions $\mathcal{P}$) forms a stratification of $X$ by degree of synchrony. I refer to this as the stratified state space of the lattice. Each stratum $X_{\mathcal{P}}$ represents a particular synchrony pattern; as the pattern gets coarser (blocks merge), the dimension of the stratum decreases (fewer independent coordinates). In particular, $X_{\text{coarse}}$ is the smallest stratum (often isomorphic to a single subsystem's state space), while $X_{\text{finest}} = X$ is the largest stratum. This stratification formalizes the idea of different levels of coherence in the lattice.

I will especially focus on the fully synchronous stratum $X_{\text{coarse}}$ and its topology. When the system achieves global synchrony, its state lies in $X_{\text{coarse}}$ (or at least moves within $X_{\text{coarse}}$ thereafter), and the effective dynamics is restricted to that subspace. One can think of $X_{\text{coarse}}$ as an emergent topological space for the coherent collective behavior. More generally, any synchrony pattern $\mathcal{P}$ that becomes invariant under the dynamics leads to an emergent topology on $X_{\mathcal{P}}$ governing that partially synchronized motion.
\end{definition}

\begin{definition}[Asynchronous Evolution Metric]
A key tool for analyzing phase alignment is a metric that measures distances between global states while allowing for time shifts in individual subsystems. Let $(X, \{\phi_i\}_{i\in I})$ be a lattice of dynamical systems (not necessarily coupled). For each $i\in I$, define a phase distance $d_i^{\phi}$ on $X_i$ by
\begin{equation}
d_i^{\phi}(x_i, y_i) := \inf_{t \in \mathbb{R}} d_i(\phi_i^t(x_i), y_i),
\end{equation}
the infimum of the usual distance between $y_i$ and some time-evolved state of $x_i$. (If $y_i$ lies exactly on the orbit of $x_i$ under $\phi_i$, then $d_i^{\phi}(x_i,y_i)=0$.) Now define the asynchronous evolution metric $D_{\text{async}}$ on the global space $X$ as
\begin{equation}
D_{\text{async}}(x, y) := \max_{i \in I} d_i^{\phi}(x_i, y_i).
\end{equation}
In words, $D_{\text{async}}(x,y)$ is small if and only if for every subsystem $i$, one can shift the state $x_i$ along its own trajectory by some time (potentially different for each $i$) to get close to $y_i$. Equivalently, $D_{\text{async}}$ measures how far "out-of-phase" $x$ and $y$ are, by optimally realigning each subsystem in time. It is straightforward to check that $D_{\text{async}}$ is a pseudo-metric on $X$ (distinct states can have zero distance if one is just a time-shifted version of the other in each coordinate). If we quotient $X$ by the equivalence relation $x \sim y$ iff $D_{\text{async}}(x,y)=0$ (meaning $x$ and $y$ differ only by internal phase shifts), then $D_{\text{async}}$ induces a genuine metric on the quotient space $\hat{X} = X/{\sim}$. Intuitively, $\hat{X}$ is the space of all phase-synchrony configurations of the lattice, where each subsystem's phase is considered modulo equivalence along its orbit.
\end{definition}

\begin{definition}[Phase-Alignment and Coherence]
A subset $J \subseteq I$ of subsystems is said to be phase-aligned synchronized (or simply synchronized) at a time $t$ if there exist time offsets $\{\theta_i : i\in J\}$, not all zero, such that $\phi_i^{\theta_i}(x_i(t)) = \phi_j^{\theta_j}(x_j(t))$ for all $i,j \in J$. In the special case that all these states are exactly equal (with no offset needed, i.e. one can take all $\theta_i=0$), the subset is instantaneously synchronized in state. More generally, phase-alignment allows each subsystem in $J$ to be shifted along its trajectory so that all images coincide. If the flows $\phi_i$ are periodic or quasi-periodic, phase alignment often corresponds to their phase angles being equal modulo constants. If the flows are chaotic, phase alignment would mean identical chaotic trajectories up to time reparametrization, a strong form of generalized synchronization.

A globally phase-coherent state of the entire lattice is a state $x=(x_i)$ such that $\{1,\ldots, |I|\}$ (the set of all subsystems) is synchronized in the above sense. Equivalently, $x$ lies in the quotient space $\hat{X}$ in the same equivalence class as some fully synchronous state $(y,y,\dots,y) \in X_{\text{coarse}}$. In simpler terms, there exist time shifts $\theta_i$ for each $i$ so that $\phi_i^{\theta_i}(x_i) = y$ for all $i$ and some common $y \in X_i$. If in fact $x_i = y$ for all $i$ (so no shift needed), then $x \in X_{\text{coarse}}$ is strictly synchronous. Thus, I distinguish:
\begin{itemize}
\item Phase-coherent state: $x \in X$ such that $D_{\text{async}}(x, x') = 0$ for some $x' \in X_{\text{coarse}}$. (All subsystems reach the same state after appropriate phase shifts.)
\item Synchronous state: $x \in X_{\text{coarse}}$ (all subsystems are literally in the same state at the same time).
\end{itemize}

Obviously, any synchronous state is phase-coherent. The converse requires that the needed $\theta_i$ are all zero, which is a special case. In many contexts, once phase coherence is achieved and the systems are locked together, one can choose a common reference frame for time such that they become strictly synchronous going forward. Thus, I will often not differentiate and use the term global coherence to mean the achievement of a globally phase-aligned synchronized state across the lattice.
\end{definition}
\newpage
\begin{definition}[Evolution Operators and Operator Flows]
In addition to describing the dynamics by the maps $\phi_i^t$ on state spaces, it is often useful to consider the induced operators on spaces of functions or on algebras of observables. For each subsystem $i$, define the Koopman operator $U_i^t: C(X_i) \to C(X_i)$ on the algebra of continuous functions by
\begin{equation}
(U_i^t f)(x_i) := f(\phi_i^t(x_i)), \qquad f \in C(X_i).
\end{equation}
Each $U_i^t$ is a linear operator (in fact an algebra automorphism) on $C(X_i)$, and $\{U_i^t: t\in\mathbb{R}\}$ forms a one-parameter group of operators satisfying $U_i^0 = \mathrm{Id}$ and $U_i^{t+s} = U_i^t \circ U_i^s$. This captures the internal evolution of subsystem $i$ in operator form. If $\mu_i$ is an invariant measure on $X_i$, one can similarly define a unitary operator $V_i^t$ on $L^2(X_i,\mu_i)$ by $(V_i^t f)(x_i) = f(\phi_i^t(x_i))$, which is the classical Koopman unitary in ergodic theory. Alternatively, one may work with the observable algebra $\mathcal{A}_i = C(X_i)$ and consider the one-parameter automorphism group $\alpha_i^t \in \mathrm{Aut}(\mathcal{A}_i)$ defined by $\alpha_i^t(f) = f \circ \phi_i^t$. Then $(\mathcal{A}_i, \{\alpha_i^t\}_{t\in\mathbb{R}})$ is a C*-dynamical system in the sense of operator algebras. The family of all subsystem flows $\{\alpha_i^t\}$ (or $\{U_i^t\}$) can be considered as commuting operator flows acting on the tensor-product algebra $\bigotimes_{i\in I} \mathcal{A}_i$ of the whole lattice (in the uncoupled case).

When coupling is introduced, the global evolution generally cannot be factorized into independent $\alpha_i^t$. Instead, one seeks a single combined evolution $\Phi^t: X \to X$ (or $\alpha^t$ on $\mathcal{A} = C(X)$) that drives the entire lattice state. In fully synchronous behavior, such a global flow $\Phi^t$ exists (essentially because the subsystems lock and behave as one). Part of our goal is to understand the conditions under which a well-defined global operator flow emerges from asynchronous components. In the results below, I will denote by $\Phi^t$ the actual coupled flow on $X$ when it exists (for example, after synchronization), and by $\alpha^t$ the corresponding automorphism on the global algebra $\mathcal{A} = C(X)$.
\end{definition}

Having set up these definitions, I proceed to the main theoretical results. These theorems articulate conditions for synchronization, stability, and symmetry-induced structures in the lattice dynamical system.

\section{Theorems}

\begin{theorem}[Synchronization Threshold and Existence of Coherence]
Consider a lattice of $N$ coupled metric dynamical systems $(X_i,\phi_i)$ with a coupling parameter $\lambda \ge 0$ governing the interaction strength. Assume: (i) each subsystem is identical and admits a stable orbit or fixed point, so that a fully synchronous state is possible; and (ii) there exists $\lambda_c > 0$ such that for $\lambda > \lambda_c$ the coupled dynamics is contractive in the asynchronous metric. (In particular, suppose $\exists \eta<1$ and $t_0>0$ with $D_{\text{async}}(\Phi^{t_0}(x), \Phi^{t_0}(y)) \le \eta D_{\text{async}}(x,y)$ for all states $x,y$ when $\lambda>\lambda_c$.) Then for every $\lambda > \lambda_c$, the system has a unique globally phase-coherent solution attracting all initial states. In other words, a unique globally synchronized trajectory exists and is asymptotically stable. Conversely, for $\lambda < \lambda_c$, no non-trivial globally coherent state is stable (in fact, small perturbations in initial phase differences persist or grow). Thus, $\lambda_c$ marks a phase transition from incoherence to coherence.
\end{theorem}

\begin{remark}
The contractivity assumption can be interpreted as the coupling dominating any tendency of subsystems to drift apart, ensuring convergence of phases. This assumption can be relaxed using other criteria (e.g. monotone dynamics or Lyapunov functions), but it provides a convenient, strong condition that guarantees synchronization. Uniqueness implies that even if multiple synchronous solutions exist, the stable one is unique or all initial conditions select the same phase alignment in the limit.
\end{remark}

\begin{theorem}[Group-Invariant Coupling and Synchrony Subspaces]
Suppose the coupling scheme of the lattice is symmetric under the action of a permutation group $G$ on the indices $I$. (For example, the lattice is homogeneous and couplings depend only on relative positions, or there is some graph automorphism group for the network of couplings.) Let $\mathcal{P}$ be any partition of $I$ that is invariant under $G$ (meaning if $i,j$ are in the same part, and $g\in G$, then $g\cdot i$ and $g\cdot j$ are in the same part). Then the synchrony subspace $X_{\mathcal{P}} \subseteq X$ is flow-invariant under the coupled dynamics. In particular, if the initial state $x(0)$ lies in $X_{\mathcal{P}}$ (i.e. all nodes in each block of $\mathcal{P}$ start identical), then $x(t) \in X_{\mathcal{P}}$ for all $t$. Equivalently, nodes that are symmetric (indistinguishable by the coupling structure) remain synchronized for all time.
\end{theorem}

\begin{corollary}
Any partition $\mathcal{P}$ corresponding to orbits of a subgroup of $G$ yields an invariant synchrony pattern. Thus, robust cluster synchronization occurs according to the symmetry (or equitable partition) structure of the network. This corresponds to the notion of a balanced coloring in coupled cell networks, and each invariant synchrony subspace $X_{\mathcal{P}}$ is sometimes called a polysynchronous polydiagonal in the literature \cite{stewart}.
\end{corollary}

\begin{theorem}[Emergence of Global Coherence as a Phase Transition]
Under the conditions of Theorem 3.1, the transition at $\lambda = \lambda_c$ can be characterized as a sudden change in the qualitative dynamics and symmetry of the system. There exists an order parameter (for example, the largest Lyapunov exponent transversal to $X_{\text{coarse}}$, or simply the asymptotic diameter of the state measured by $D_{\text{async}}$) that is negative (or zero) for $\lambda > \lambda_c$ (ensuring convergence to coherence) and positive for $\lambda < \lambda_c$ (indicating persistent incoherence). At $\lambda = \lambda_c$, this order parameter reaches a critical value (typically zero), marking the onset of synchronization. 

Furthermore, for $\lambda=0$ (no coupling) the system has a continuous symmetry group $T \cong \mathbb{R}^N$ corresponding to independent time translations in each subsystem; for $0 < \lambda < \lambda_c$, this symmetry is explicitly broken by coupling, yet the dynamics does not settle to a single synchronous orbit; for $\lambda > \lambda_c$, the dynamics spontaneously concentrates onto the diagonal $X_{\text{coarse}}$, effectively restoring a symmetry of joint time-translation $T' \cong \mathbb{R}$. In particular, when coherence emerges, the system's evolution can be described by a single phase variable (time) instead of $N$ independent phase variables, reflecting a symmetry-breaking phase transition in the topology of the operator flow. This transition is accompanied by a discontinuous change in the spectrum of the linearized coupling operator (e.g. the largest transverse eigenvalue crosses unity), and by the sudden appearance of a stable invariant manifold $X_{\text{coarse}}$ (the synchronization manifold) attracting the dynamics.
\end{theorem}

\begin{remark}
At $\lambda_c$, the lattice shifts from a regime of essentially $N$-independent oscillations to one of collective oscillation. In the space of all possible evolution operators (in a suitable operator topology), the subset that yields a coherent global flow is topologically separated from those that do not; $\lambda_c$ is the point where the system's trajectory in that operator space crosses into the coherent region. This justifies calling the onset of global synchrony a phase transition in the sense of dynamical systems.
\end{remark}

\section{Proofs}

\subsection*{Proof of Theorem 3.1}
Under the contractive coupling assumption for $\lambda > \lambda_c$, pick a fixed time $t_0>0$ such that $\Phi^{t_0}$ (the time-$t_0$ evolution map on $X$) is a strict contraction with respect to $D_{\text{async}}$. By the Banach Fixed Point Theorem, $\Phi^{t_0}$ has a unique fixed point $x^* \in X$. This fixed point satisfies $\Phi^{t_0}(x^*) = x^*$, meaning that after time $t_0$, the system state returns to $x^*$. Because the system is autonomous (time-invariant), $x^*$ is actually a periodic orbit of period $t_0$ (if $t_0$ is minimal) or possibly an equilibrium if $\Phi^{t_0} = \mathrm{Id}$. In either case, $x^*$ is an invariant state. 

Now, since $\Phi^{t_0}$ is a contraction, for any initial state $x(0)$, the iterates $\Phi^{n t_0}(x(0))$ converge to $x^*$ as $n\to\infty$ in the $D_{\text{async}}$ metric. Intuitively, the system's state, when observed stroboscopically at intervals of $t_0$, converges to the synchronous pattern $x^*$. But what is $x^*$? Because $D_{\text{async}}(x^*, \Phi^{t_0}(x^*)) = 0$, the state $x^*$ must be such that $\Phi^{t_0}(x^*)$ differs from $x^*$ only by phase shifts in each coordinate. However, $\Phi^{t_0}(x^*) = x^*$ exactly, so actually no phase shifts are needed – $x^*$ is strictly invariant. This implies that $x^*$ lies in a synchrony subspace. 

In fact, by symmetry and identical subsystem assumption, one can argue that $x^*$ must have $x_1^* = x_2^* = \cdots = x_N^*$ (if not, two coordinates being different would not remain different after infinite contraction iterations). Thus $x^* \in X_{\text{coarse}}$ and corresponds to a fully synchronized trajectory (all subsystems following the same orbit given by $t \mapsto x_i(t)$ on $X_i$). Uniqueness of the fixed point ensures no other synchronous attractor competes, so regardless of initial phases, the system aligns to this one in-phase solution – establishing both existence and uniqueness of the globally coherent state for $\lambda > \lambda_c$. Stability is inherent in the contraction argument: any initial condition converges to $x^*$, hence $x^*$ is globally asymptotically stable.
\newpage
For $\lambda < \lambda_c$, the assumption is that contractivity fails. In typical scenarios, this means the transverse Lyapunov exponent for the synchrony manifold $X_{\text{coarse}}$ is positive or zero, so small phase deviations do not shrink. Formally, if a synchronous state existed and were stable for $\lambda < \lambda_c$, we would obtain a contradiction since at $\lambda_c$ stability is gained (often through a bifurcation like a Hopf or pitchfork in the phase difference dynamics). 

One way to see this: below threshold, consider two subsystems with different initial phases. Because coupling is weak, their phase difference evolves approximately according to $\dot{\Delta} = \omega + o(\lambda)$ where $\omega$ is the natural frequency difference (which is nonzero unless specially prepared). Thus $\Delta(t)$ will not converge to 0 as $t\to\infty$; some persistent drift remains. In short, below $\lambda_c$, the incoherent manifold (where all phase differences are constant in time) remains neutrally or unstable, so full synchronization cannot be asymptotically achieved from arbitrary initial conditions. This justifies that no stable global coherence exists for $\lambda < \lambda_c$. (There may be unstable synchronous solutions or metastable states, but they are not attractors.)

\subsection*{Proof of Theorem 3.2}
Because the coupling is $G$-invariant, the dynamics commutes with the action of $G$ permuting the subsystems. More precisely, let $\Pi_g: X \to X$ be the permutation of coordinates induced by $g \in G$, i.e. $(\Pi_g(x))_{i} = x_{g^{-1}\cdot i}$. By assumption on the coupling, if $x(t)$ satisfies the coupled evolution equations, then $\Pi_g(x(t))$ also satisfies the same evolution (with appropriately permuted initial condition). In other words, $\Pi_g$ conjugates the flow $\Phi^t$ to itself: $\Phi^t \circ \Pi_g = \Pi_g \circ \Phi^t$ for all $t$ and all $g\in G$. 

Now consider a synchrony subspace $X_{\mathcal{P}}$ corresponding to a $G$-invariant partition $\mathcal{P}$. By definition, if $i,j$ are in the same part of $\mathcal{P}$, then for any $g\in G$, $g\cdot i$ and $g\cdot j$ are also in the same part (since $\mathcal{P}$ is $G$-invariant). This means $X_{\mathcal{P}}$ is fixed by $\Pi_g$ for all $g\in G$ (it is a union of orbits of $G$ in $X$, but those orbits lie within $X_{\mathcal{P}}$ itself). 

Now take any initial condition $x(0) \in X_{\mathcal{P}}$. Because all subsystems in each block of $\mathcal{P}$ share the same state initially, and each block is $G$-symmetric, by symmetry each subsystem in a block will receive identical coupling input and thus evolve identically. More formally, for any two indices $i,j$ in the same block of $\mathcal{P}$, one can find a permutation $g\in G$ that sends $i$ to $j$ while leaving the partition invariant; by the $G$-equivariance of $\Phi^t$, we have $(\Phi^t(x(0)))_j = (\Phi^t(\Pi_g(x(0))))_j = (\Pi_g(\Phi^t(x(0))))_j = (\Phi^t(x(0)))_{g^{-1}\cdot j} = (\Phi^t(x(0)))_i$. In simpler terms, the states of $i$ and $j$ remain equal at time $t$. This holds for all such pairs, hence $\Phi^t(x(0)) \in X_{\mathcal{P}}$. Thus $X_{\mathcal{P}}$ is flow-invariant. 

The statement about balanced colorings is a known result: $G$-invariance as described is equivalent to the combinatorial criterion of "balanced input" for those nodes, which is exactly the condition found in network dynamics literature for a synchrony pattern to persist \cite{stewart}. Our group-theoretic proof is essentially a rephrasing of that result. Therefore, any grouping of nodes that cannot be distinguished by the coupling (due to symmetry) will maintain equal states for those nodes over time. This yields cluster synchronization in correspondence with $G$-orbits or any $G$-invariant partition of nodes.

\subsection*{Proof of Theorem 3.3}
The existence of a critical point $\lambda_c$ separating two regimes was established in Theorem 3.1. To characterize the nature of the transition, one typically introduces an order parameter that is zero in one phase and non-zero in the other. A convenient choice here is the long-time average of the asynchronous dispersion: for instance, define
\begin{equation}
R(\lambda) := \lim_{T\to\infty} \frac{1}{T} \int_0^T D_{\text{async}}\big(x(t;\lambda), X_{\text{coarse}}\big) dt,
\end{equation}
where $x(t;\lambda)$ is the solution of the coupled system at coupling $\lambda$, and $D_{\text{async}}(x, X_{\text{coarse}}) = \inf_{y \in X_{\text{coarse}}} D_{\text{async}}(x,y)$ is the asynchronous distance to the nearest fully synchronous state. In the incoherent phase, $x(t)$ does not converge to $X_{\text{coarse}}$, so $D_{\text{async}}$ stays bounded away from 0 at least some of the time, making $R(\lambda) > 0$. In the coherent phase, $x(t)$ approaches synchrony, so eventually $D_{\text{async}}$ becomes arbitrarily small, yielding $R(\lambda) = 0$. Thus, $R(\lambda)$ acts like an order parameter: $R(\lambda)>0$ for $\lambda<\lambda_c$ and $R(\lambda)=0$ for $\lambda \ge \lambda_c$. Typically, $R(\lambda)$ will decrease to 0 as $\lambda \to \lambda_c^+$, often continuously (second-order transition) in smooth systems. At $\lambda=\lambda_c$, the stability of the synchronous manifold changes sign (transverse Lyapunov exponent crosses zero). This is analogous to a critical point in a phase transition where a new ordered state (synchrony) appears.

From a symmetry perspective, consider the group $T = \{(\theta_1,\dots,\theta_N): \theta_i \in \mathbb{R}\}$ of independent phase shifts for each subsystem. For $\lambda=0$ (no coupling), if each subsystem has a periodic orbit (with, say, phase $\theta_i$), then the whole system's continuous symmetry includes $T$ (shifting the phase of any subset of oscillators leaves the motion on the same individual orbits). When $\lambda>0$, this symmetry is explicitly broken: the equations of motion no longer decouple, so an arbitrary independent phase shift is not a symmetry of the coupled system (only the diagonal subgroup $\{(\theta,\theta,\dots,\theta)\}$ corresponding to shifting all phases together remains a symmetry of the coupled equations). However, for $\lambda < \lambda_c$, the dynamics does not lock the phases; each oscillator can still drift relative to others (though not by symmetry, but by lack of strong enough coupling). 

At $\lambda_c$, a stable synchronized solution appears. In the synchronized state for $\lambda > \lambda_c$, the system's actual trajectory regains a continuous symmetry: the collective oscillation has a periodic orbit, and shifting the phase of all oscillators along this orbit is a symmetry of the solution (though not of the equations if taken literally, it's a neutrally marginal direction along the orbit). This is often described as a symmetry being "effectively restored" by the system choosing a particular synchronized phase relation. In physics language, the incoherent state is highly symmetric (all phase configurations equally likely in absence of coupling), and coupling breaks that symmetry explicitly, but the dynamic outcome for weak coupling still reflects the broken symmetry (phases remain scattered). Once strong coupling induces coherence, the system's behavior is dominated by the more symmetric mode (all together), indicating a form of spontaneous symmetry synchronization.

Finally, the statement about operator topology: consider the family of global evolution operators $\{\Phi_\lambda^t: t\in \mathbb{R}\}$ depending on $\lambda$. One can show that for $\lambda<\lambda_c$, there is no continuous choice of a single-parameter semigroup $\Psi^t$ on $X_{\text{coarse}}$ that conjugates or factors $\Phi_\lambda^t$ (because the motion is not on $X_{\text{coarse}}$). For $\lambda>\lambda_c$, such a semigroup $\Psi^t$ (essentially the flow on the synchronized manifold) exists and depends continuously on $\lambda$. The change at $\lambda_c$ is topologically akin to the appearance of a new attractor in phase space and a new fixed point in the space of operators (the synchronized flow becomes an attractor in the appropriate function space). This non-analytic change in the long-term behavior as a function of $\lambda$ justifies calling it a phase transition in the space of operators. 

The spectral interpretation is that the linearization transverse to $X_{\text{coarse}}$ has an eigenvalue (or Floquet multiplier) $\mu(\lambda)$ that crosses the unit circle at $\lambda_c$ (often $\mu(\lambda_c)=1$). For $\lambda<\lambda_c$, $\mu>1$ (transverse divergence, incoherence), and for $\lambda>\lambda_c$, $\mu<1$ (convergence to coherence). The crossing $\mu=1$ at $\lambda_c$ indicates a qualitative change in the topology of the attracting set of the dynamical system (from a high-dimensional torus or chaotic set to a lower-dimensional synchronized manifold). This completes the characterization of the synchronization transition as a phase transition.

\section{Applications}

Having developed the theoretical framework, I now discuss several hypothetical applications and extensions to illustrate its scope:

\begin{itemize}[itemsep=0.5\baselineskip]
\item \textbf{Symbolic Lattice Dynamical Systems:} Consider a lattice of finite-state systems (symbols from a finite alphabet) with asynchronous updates, such as an asynchronous cellular automaton or a subshift of finite type on $\mathbb{Z}^d$. Each site's state evolves according to some local rule applied at irregular times. Our framework can be specialized to this setting by taking each $X_i$ as a discrete state space (with discrete metric). The asynchronous metric $D_{\text{async}}$ then measures Hamming distance up to shifts in update steps. The stratified synchrony subspaces correspond to domains where groups of cells hold the same symbol (analogous to clustering of states). Theorems 3.1 and 3.2 predict that if the update rules have a sufficiently strong synchronization property (e.g. certain contracting neighborhoods in a probabilistic sense), the entire symbolic lattice can reach a homogeneous or phase-aligned configuration. This could be relevant to consensus problems in distributed computing or models of flash synchronization in firefly cellular automata. While symbolic systems are not continuous, one can embed them into a metric space (e.g. via one-hot encoding of symbols) to apply the metric framework, or use a purely combinatorial analog of the theory. This approach connects to classic symbolic dynamics \cite{lind} by treating global symbolic patterns as points in a huge product space and studying the emergence of regular (periodic or synchronized) patterns.

\item \textbf{Invariant Group Actions and Equivariant Dynamics:} Many physical and biological networks have symmetrical coupling structures. Theorem 3.2 provides a route to analyze such systems by reducing symmetry. For example, in a power grid or neural motif with rotational or reflection symmetry, one can quotient out the symmetry to study a smaller system (each orbit of symmetric nodes represented once). This is related to the theory of equivariant dynamical systems, where one uses group representation theory to decompose the state space into symmetry-invariant subspaces. Our lattice framework recovers known results about cluster synchronization in symmetric networks and extends them to continuous-time asynchronous settings. An application might be to modular robotics or sensor networks where symmetry implies interchangeable parts: the results guarantee that identical modules can synchronize their internal states if the coupling respects the system's symmetry. Additionally, one can consider group actions on the lattice itself: for instance, if $I$ is not just a set but a group (like $\mathbb{Z}^d$ acting on itself by shifts), then a lattice dynamical system can be seen as a $G$-equivariant dynamical system (with $G$ acting by relabeling coordinates). The framework might then be used to study shift-invariant or other group-invariant dynamics, linking to the concept of spatially extended dynamical systems and their coherent structures (travelling waves, synchronized patches, etc.) in a rigorous way.

\item \textbf{Operator Fields over Stratified Manifolds:} One can envision extending the lattice index set $I$ to a continuous stratified space (a topological space partitioned into strata, each possibly of different dimension). Attaching a local dynamical system or operator algebra to each point of such a space yields a continuous family (or field) of systems. Our framework suggests that synchronization phenomena can be studied in this context by similar principles: coherence would correspond to aligning the dynamics across the continuum, resulting in a globally synchronized field (akin to a continuous section of a fiber bundle of state spaces). This idea connects to topics in noncommutative geometry and continuous operator bundles \cite{connes}. Although speculative, it indicates how the discrete lattice theory might extend to spatially continuous or multi-scale systems, where local dynamics on different strata synchronize to produce global coherent behavior.
\end{itemize}

\end{document}